\documentclass[a4paper, 11pt, reqno]{amsart}

\usepackage[margin=25truemm]{geometry}
\usepackage{hyperref}

\title{Quantum $K$-theoretic Whitney relation for Lagrangian Grassmannian}
\author{Takafumi Kouno}
\address{Department of Applied Mathematics, Faculty of Science and Engineering, Waseda University, 3-4-1 Okubo, Shinjuku-ku, Tokyo 169-8555, Japan}
\email{t.kouno@aoni.waseda.jp}

\numberwithin{equation}{section}

\usepackage{amssymb}
\usepackage{mathtools}
\mathtoolsset{showonlyrefs}
\usepackage{enumitem}
\usepackage{dynkin-diagrams}

\newcommand{\BZ}{\mathbb{Z}}
\newcommand{\BR}{\mathbb{R}}
\newcommand{\BC}{\mathbb{C}}
\newcommand{\udBC}{\underline{\BC}}
\newcommand{\CA}{\mathcal{A}}
\newcommand{\CE}{\mathcal{E}}
\newcommand{\CF}{\mathcal{F}}
\newcommand{\CO}{\mathcal{O}}
\newcommand{\CP}{\mathcal{P}}
\newcommand{\CQ}{\mathcal{Q}}
\newcommand{\CS}{\mathcal{S}}
\newcommand{\FS}{\mathfrak{S}}
\newcommand{\Fh}{\mathfrak{h}}
\newcommand{\ve}{\varepsilon}
\newcommand{\vpi}{\varpi}
\newcommand{\e}{\mathbf{e}}
\newcommand{\vtl}{\vartriangleleft}
\newcommand{\cl}{\mathsf{cl}}
\newcommand{\Sp}[1]{\operatorname{Sp}_{#1}(\BC)}
\newcommand{\LG}{\operatorname{LG}}
\newcommand{\Fl}{\operatorname{Fl}}
\newcommand{\FlC}[1]{\operatorname{Fl}^{C_{#1}}}
\newcommand{\bra}[1]{[\![ #1 ]\!]}
\newcommand{\class}[1]{\left[ #1 \right]}
\newcommand{\mcr}[1]{\lfloor #1 \rfloor}
\newcommand{\wh}[1]{\widehat{#1}}
\newcommand{\ti}[1]{\widetilde{#1}}
\newcommand{\pair}[2]{\langle #1, #2 \rangle}

\newcommand{\extprod}{\bigwedge\nolimits}

\DeclareMathOperator{\rk}{rk}
\DeclareMathOperator{\ed}{end}
\DeclareMathOperator{\down}{down}
\DeclareMathOperator{\wt}{wt}
\DeclareMathOperator{\QBG}{QBG}
\DeclareMathOperator{\Lie}{Lie}
\DeclareMathOperator{\diag}{diag}

\theoremstyle{plain}
\newtheorem{thm}{Theorem}[section]
\newtheorem{lem}[thm]{Lemma}
\newtheorem{prop}[thm]{Proposition}
\newtheorem{cor}[thm]{Corollary}

\theoremstyle{definition}
\newtheorem{defn}[thm]{Definition}

\theoremstyle{remark}
\newtheorem{rem}[thm]{Remark}
\newtheorem{exm}[thm]{Example}

\newenvironment{enu}{%
\begin{enumerate}[label=\textup{(\arabic*)}, left=0pt]}{%
\end{enumerate}}

\begin{document}

\begin{abstract}
We give a quantum $K$-theoretic analogue of the Whitney relation for the torus-equivariant $K$-ring of the Lagrangian Grassmannian. As a consequence, we give a presentation of the torus-equivariant quantum $K$-ring of the Lagrangian Grassmannian. The strategy for obtaining such relations is to use the quantum $K$-theoretic Whitney relation for the full-flag manifold in type $C$ established by the author and Kato's pushforward map from the quantum $K$-ring of the full-flag manifold to that of the partial flag manifold. Then the Nakayama-type argument established by Gu et al., together with the classical Whitney presentation for the torus-equivariant $K$-ring of the Lagrangian Grassmannian, gives a proof of such a presentation. 
\end{abstract}

\maketitle

\section{Introduction}

Let $\Sp{2n}$ be the symplectic group with respect to a fixed non-degenerate skew-symmetric form $(-, -)$ on $\BC^{2n}$, and take a maximal torus $T \subset \Sp{2n}$. 
We denote by $R(T)$ the representation ring of $T$. 
The \emph{Lagrangian Grassmannian} $\LG(n)$ is defined as 
\begin{equation}
\LG(n) := \left\{ F \ \middle| \ \text{$F$ is an isotropic subspace of $\BC^{2n}$ with respect to $(-, -)$ with $\dim_{\BC}F = n$} \right\}. 
\end{equation}
Our interest is to give a presentation of the \emph{$T$-equivariant quantum $K$-ring} $QK_{T}(\LG(n))$ of the Lagrangian Grassmannian, written as 
\begin{equation}
QK_{T}(\LG(n)) = K_{T}(\LG(n)) \otimes_{R(T)} R(T)\bra{Q}. 
\end{equation}
For this purpose, we focus on the Whitney presentation. 
For a $T$-equivariant vector bundle $\CE$ on $\LG(n)$, we define the \emph{Hirzebruch's $\lambda_{y}$-classes} $\lambda_{y}(\CE) \in K_{T}(\LG(n))[y]$ by 
\begin{equation}
\lambda_{y}(\CE) := \sum_{d=0}^{\rk(\CE)} y^{d} \class{\extprod^{d} \CE}, 
\end{equation}
where $\rk(\CE)$ denotes the rank of $\CE$. 
Then it is well-known that for an exact sequence $0 \rightarrow \CE_{1} \rightarrow \CE_{2} \rightarrow \CE_{3} \rightarrow 0$ of $T$-equivariant vector bundles on $\LG(n)$, we have $\lambda_{y}(\CE_{2}) = \lambda_{y}(\CE_{1}) \cdot \lambda_{y}(\CE_{3})$. 
This relation is called a \emph{Whitney relation}. 
Hence we obtain 
\begin{equation} \label{eq:classical_Whitney}
\lambda_{y}(\CS) \cdot \lambda_{y}(\CQ) = \lambda_{y}(\udBC^{2n}),  
\end{equation}
where $\CS$ (resp., $\CQ$) denotes the tautological bundle (resp., the universal quotient bundle) on $\LG(n)$, while $\udBC^{2n}$ denotes the trivial bundle on $\LG(n)$ of rank $2n$. 
Moreover, it is well-known that \eqref{eq:classical_Whitney} gives the complete list of defining relations of the $T$-equivariant (classical) $K$-ring $K_{T}(\LG(n))$. 
This presentation is called the \emph{Whitney presentation} of $K_{T}(\LG(n))$. 

One of the main results of this paper is to give a quantum analogue of \eqref{eq:classical_Whitney}, which we call the \emph{quantum $K$-theoretic Whitney relation} (or the \emph{$QK$-Whitney relation} for short).
Note that this relation was conjectured by Gu, Mihalcea, Sharpe, and Zou (\cite{GMSZ}) from a physics perspective. 
\begin{thm}[{Conjectured in \cite[(6.95)]{GMSZ}}] \label{thm:Whitney_rel}
In $QK_{T}(\LG(n))[y]$, we have 
\begin{equation} \label{eq:QK-Whitney_rel}
\lambda_{y}(\CS) \star \lambda_{y}(\CQ) = \lambda_{y}(\udBC^{2n}) - y^{2} Q \sum_{d=0}^{2n-2} y^{d} \CF_{d}, 
\end{equation}
where $\CF_{d} \in K_{T}(\LG(n))$ is defined as 
\begin{equation}
\CF_{d} := \class{\extprod^{d} \CQ} + \class{\extprod^{d-2} \CQ} + \cdots + \class{\extprod^{d - 2\lfloor d/2 \rfloor} \CQ}. 
\end{equation}
for $0 \le d \le n-1$, and $\CF_{d} := \CF_{2n-2-d}$ for $n \le d \le 2n-2$. 
\end{thm}

As a consequence, we obtain a presentation of $QK_{T}(\LG(n))$. 
Let $e_{d}$ denote the $d$-th elementary symmetric polynomial for $d \in \mathbb{Z}_{\ge 0}$. 
Set $e_{d}(x) := e_{d}(x_{1}, \ldots, x_{n})$, and $e_{d}(y) := e_{d}(y_{1}, \ldots, y_{n})$ for $0 \le d \le n$. 
Let $\{\ve_{k} \mid 1 \le k \le n\}$ be the standard basis of the weight lattice $\Lambda$ of $\Sp{2n}$. 
We denote by $\e^{\lambda}$ the character of $T$ associated to $\lambda \in \Lambda$, and set 
\begin{equation}
e_{d}(\e^{\pm \ve}) := e_{d}(\e^{\ve_{1}}, \ldots, \e^{\ve_{n}}, \e^{-\ve_{n}}, \ldots, \e^{-\ve_{1}})
\end{equation}
for $0 \le d \le 2n$. 
Let $\FS_{n}$ be the $n$-th symmetric group. Then $\FS_{n}$ acts on the polynomial ring $(R(T)\bra{Q})[x_{1}, \ldots, x_{n}, y_{1}, \ldots, y_{n}]$ as the permutation of indices of variables. 
Denote by $(R(T)\bra{Q})[x_{1}, \ldots, x_{n}, y_{1}, \ldots, y_{n}]^{\FS_{n}}$ the subring of $(R(T)\bra{Q})[x_{1}, \ldots, x_{n}, y_{1}, \ldots, y_{n}]$ composed of all $\FS_{n}$-invariants; i.e., 
\begin{equation}
(R(T)\bra{Q})[x_{1}, \ldots, x_{n}, y_{1}, \ldots, y_{n}]^{\FS_{n}} = (R(T)\bra{Q})[e_{d}(x), e_{d}(y) \mid 1 \le d \le n]. 
\end{equation} 
Define $J_{\LG(n)}$ as an ideal of $(R(T)\bra{Q})[x_{1}, \ldots, x_{n}, y_{1}, \ldots, y_{n}]^{\FS_{n}}$ generated by 
\begin{equation}
\sum_{\substack{k, l \ge 0 \\ k+l=d}} e_{k}(x) e_{l}(y) - e_{d}(\e^{\pm \ve}) + QF_{d-2}, \qquad 1 \le d \le 2n, 
\end{equation}
where $F_{d}$ is defined as 
\begin{equation}
F_{d} := \sum_{k=0}^{\lfloor d/2 \rfloor} e_{d-2k}(y)
\end{equation}
for $0 \le d \le n-1$, $F_{d} := F_{2n-2-d}$ for $n \le d \le 2n-2$, and $F_{-1} := 0$. 

\begin{thm} \label{thm:Whitney_presentation}
We have an $R(T)\bra{Q}$-algebra isomorphism
\begin{equation}
\Psi: (R(T)\bra{Q})[x_{1}, \ldots, x_{n}, y_{1}, \ldots, y_{n}]^{\mathfrak{S}_{n}}/J_{\LG(n)} \xrightarrow{\sim} QK_{T}(\LG(n))
\end{equation}
such that $\Psi(e_{d}(x)) = \class{\extprod^{d} \CS}$ and $\Psi(e_{d}(y)) = \class{\extprod^{d} \CQ}$ for $1 \le d \le n$. 
\end{thm}

The main idea of the proof of Theorem~\ref{thm:Whitney_rel} is to utilize the pushforward $QK_{T}(\Sp{2n}/B) \twoheadrightarrow QK_{T}(\LG(n))$ of $R(T)$-algebras, established by Kato \cite{Kato}; where $B$ denotes a Borel subgroup of $G$ such that $T \subset B$. 
In \cite{Kouno}, the author obtained the $QK$-Whitney relations for the full-flag manifold $\Sp{2n}/B$ in type $C$. 
By applying Kato's pushforward to such relations, we obtain the desired relation in $QK_{T}(\LG(n))$. 
However, note that to use Kato's pushforward, we must expand elements in $QK_{T}(\Sp{2n}/B)$ as a linear combination of Schubert classes. 
For this purpose, we use the Chevalley formula computed by Lenart, Naito, and Sagaki \cite{LNS}, and we also use the classical $K_{T}$-Chevalley formula in $K_{T}(G/B)$ computed by Lenart and Postnikov \cite{LP}. 
Note that we should pay attention not to cause indeterminate forms of the form $0/0$ in the computation of Kato's pushforward, because Kato's pushforward requires substitution $Q_{k} = 1$ for $1 \le k \le n-1$ and usually $\class{\extprod^{d} \CS}$ and $\class{\extprod^{d} \CQ}$ have factors of the form $1/(1-Q_{k})$.

We emphasize that images of classes of line bundles via Kato's pushforward are non-trivial in general. 
In our Lagrangian Grassmannian case, we can describe that such images vanish (Lemma~\ref{lem:line_proj}), or be written as single classes of exterior products of $\CQ$ (Corollary~\ref{cor:quotient}). 
This enables us to find key recursive identities \eqref{eq:recursive} in $QK_{T}(\LG(n))[y]$, and finally we obtain the main relation \eqref{eq:QK-Whitney_rel}.

To prove Theorem~\ref{thm:Whitney_presentation}, we utilize the Nakayama-type argument established by Gu, Mihalcea, Sharpe, Xu, Zhang, and Zou \cite{GMSXZZ2}. 
This argument enables us to extend an $R(T)$-algebra isomorphism whose image is $K_{T}(\LG(n))$ to some $R(T)\bra{Q}$-algebra isomorphism whose image is $QK_{T}(\LG(n))$ if we have ``enough'' relations in $QK_{T}(\LG(n))$. 
Therefore, if $K_{T}(\LG(n))$ is equipped with a presentation, then we can obtain a presentation of $QK_{T}(\LG(n))$. 
Note that Nakayama's lemma inspires this argument in the theory of commutative rings. 
Again, it is known that the Whitney relation \eqref{eq:classical_Whitney} gives a presentation of $K_{T}(\LG(n))$. 
By applying this argument to \eqref{eq:QK-Whitney_rel}, we obtain the desired presentation of $QK_{T}(\LG(n))$. 

Now, we would like to mention previous studies related to the $QK$-Whitney presentation in type $A$. 
Let $\Fl^{A_{n}}(r_{1}, \ldots, r_{l})$ be the partial flag manifold in type $A_{n}$ for a certain sequence $1 \le r_{1} < \cdots < r_{l} \le n$ of integers. 
Then we can take a tautological filtration $0 = \CS_{0} \subset \CS_{1} \subset \cdots \subset \CS_{l} \subset \CS_{l+1} = \udBC^{n+1}$ of the trivial vector bundle $\udBC^{n+1}$ over $\Fl^{A_{n}}(r_{1}, \ldots, r_{l})$. Huq-Kuruvilla \cite{HK} computed the product $\lambda_{y}(\CS_{k}) \star \lambda_{y}(\CS_{k+1}/\CS_{k})$ for $k = 1, 2, \ldots, l$; note that this was first conjectured by Gu, Mihalcea, Sharpe, Xu, Zhang, and Zou \cite{GMSXZZ}. By Gu, Mihalcea, Sharpe, Xu, Zhang, and Zou \cite{GMSXZZ2}, it was proved that these products provide defining relations of $QK_{T}(\Fl^{A_{n}}(r_{1}, \ldots, r_{l}))$, by utilizing the Nakayama-type argument. 
Note that the argument to show Theorem~\ref{thm:Whitney_rel} in this paper is completely different from that used in \cite{HK}. On the other hand, we use the same argument as \cite{GMSXZZ2} to show Theorem~\ref{thm:Whitney_presentation}. 

This paper is organized as follows. 
In Section~\ref{sec:preliminary}, we fix basic notation and then recall the $QK$-Whitney relation for the flag manifold of type $C$. 
In Section~\ref{sec:Kato_push}, we review Kato's pushforward. Then, by utilizing quantum alcove models and Chevalley formulas, we compute images of line bundle classes via Kato's pushforward. 
In Section~\ref{sec:recursive}, we give a proof of Theorem~\ref{thm:Whitney_rel} by providing a certain recursive identity in $QK_{T}(\LG(n))[y]$. 
In Section~\ref{sec:presentation}, we give a proof of Theorem~\ref{thm:Whitney_presentation} by utilizing the Nakayama-type argument.

\subsection*{Acknowledgement}
The author would like to thank Takeshi Ikeda and Leonardo Mihalcea for helpful suggestions. 
The author was supported by JSPS Grant-in-Aid for Early-Career Scientists 26K16968. 

\section{Preliminary} \label{sec:preliminary}
In this section, we fix basic settings and notation. Then we recall the $QK$-Whitney relation for the flag manifold of type $C$.  
\subsection{Notation of algebraic groups and root systems}
Let $G$ be a connected, simply-connected, simple algebraic group over $\BC$ of rank $n$. 
Take a maximal torus $T \subset G$ and a Borel subgroup $B \subset G$ containing $T$, 
and take an opposite Borel subgroup $B^{-} \subset G$ such that $B \cap B^{-} = T$. 
Let $\Fh$ be the Lie algebra of $T$. we denote by $\pair{\cdot}{\cdot}: \Fh^{\ast} \times \Fh \rightarrow \BC$ the canonical pairing, where $\Fh^{\ast} = \operatorname{Hom}_{\BC}(\Fh, \BC)$. 

Take a root system $\Delta$ of $G$, and the set $\Delta^{+}$ of positive roots. 
Let $I = \{1, \ldots, n\}$ be the set of nodes of the Dynkin diagram of $G$. 
For $i \in I$, we denote by $\alpha_{i}$ the $i$-th simple root. 
For $\alpha \in \Delta$, the corresponding coroot is written as $\alpha^{\vee}$. 
Then the coroot lattice $Q^{\vee}$ of $G$ is defined as $Q^{\vee} := \bigoplus_{i \in I} \BZ \alpha_{i}^{\vee}$, where its positive part $Q^{\vee, +}$ is defined as $Q^{\vee, +} := \sum_{i \in I} \BZ_{\ge 0} \alpha_{i}^{\vee}$. 
For $\alpha \in \Delta$, we set 
\begin{equation}
|\alpha| := \begin{cases}
\alpha & \text{if $\alpha \in \Delta^{+}$}, \\ 
-\alpha & \text{if $\alpha \in -\Delta^{+}$}. 
\end{cases}
\end{equation}

Let $\Lambda = \bigoplus_{i \in I} \BZ \vpi_{i}$ be the weight lattice of $G$, where $\vpi_{i}$ denotes the $i$-th fundamental weight. 
For $\lambda \in \Lambda$, we denote by $\BC_{\lambda}$ the one-dimensional representation of $B$ whose weight is $\lambda$, and by $\e^{\lambda}$ the character of $\BC_{\lambda}$. 
The representation ring $R(T)$ of $T$ is written as $R(T) = \BZ[\e^{\lambda} \mid \lambda \in \Lambda]$. 

For the subset $L \subset I$, we take the standard parabolic subgroup $(B \subset) \, P_{L} \subset G$ corresponding to $L$. 
Then a quotient $G/P_{L}$ is called a \emph{partial flag manifold}; if $P_{L} = B$, i.e., $L = \emptyset$, a quotient $G/B$ is called the \emph{flag manifold}. 

Let $W = \langle s_{i} \mid i \in I \rangle$ be the Weyl group of $G$, where $s_{i}$, $i \in I$, denotes the simple reflection corresponding to $\alpha_{i}$. 
For $w \in W$, set $X^{w} := \overline{B^{-}wB/B}$, called the (opposite) \emph{Schubert subvariety} of $G/B$. 
Its structure sheaf is denoted by $\CO^{w} := \CO_{X^{w}}$ for $w \in W$. 
Also, we denote by $\CO(\lambda)$ the line bundle $G \times^{B} \BC_{-\lambda} \twoheadrightarrow G/B$ corresponding to $\lambda \in P$. 
In addition, for $L \subset I$, define $W_{L} := \langle s_{i} \mid i \in L \rangle \subset W$ and $W^{L} := \{ w \in W \mid \text{$\ell(ws_{i}) > \ell(w)$ for $i \in L$} \} \simeq W/W_{L}$, where $\ell(\cdot)$ denotes the length function on $W$. The \emph{Schubert subvariety} of $G/P_{L}$ is defined as $X_{L}^{w} := \overline{B^{-}wP_{L}/P_{L}}$ for $w \in W^{L}$. Its structure sheaf is denoted by $\CO_{L}^{w} := \CO_{X_{L}^{w}}$ for $w \in W^{L}$. 

\subsection{\texorpdfstring{Quantum $K$-ring}{Quantum K-ring}}
We recall the notion of the torus-equivariant quantum $K$-ring of a partial flag manifold introduced by Givental \cite{Givental} and Lee \cite{Lee}. 
Consider the case $L = I \setminus \{i_{1} < \cdots < i_{k}\}$ and set $X = G/P_{L}$. 
Let $K_{T}(X)$ be the $T$-equivariant $K$-ring of $X$. 
Note that since $X$ is smooth, $K_{T}(X)$ coincides with an $R(T)$-algebra generated by formal classes $[\CE]$ for $T$-equivariant vector bundles $\CE$ on $X$ with relations $[\CE_{2}] = [\CE_{1}] + [\CE_{3}]$ for all exact sequences $0 \rightarrow \CE_{1} \rightarrow \CE_{2} \rightarrow \CE_{3} \rightarrow 0$ of $T$-equivariant vector bundles on $X$. 
The \emph{$T$-equivariant quantum $K$-ring} of $X$ is, as an $R(T)$-module, defined by 
\begin{equation}
QK_{T}(X) := K_{T}(X) \otimes_{R(T)} R(T)\bra{Q_{i_{1}}, \ldots, Q_{i_{k}}}, 
\end{equation}
where $Q_{i_{1}}, \ldots, Q_{i_{k}}$ are (Novikov) variables. For $\xi = \sum_{1 \le j \le k} c_{j} \alpha_{i_{j}}^{\vee} \in Q^{\vee, +}_{L} := \sum_{1 \le j \le k} \BZ_{\ge 0} \alpha_{i_{j}}^{\vee}$, write $Q^{\xi} := \prod_{1 \le j \le k} Q_{i_{j}}^{c_{j}}$. 
Note that the product $\star$ of $QK_{T}(X)$ is given by \emph{$K$-theoretic Gromov--Witten invariants}, of which we omit an explanation. 
It is well-known that the set of \emph{Schubert classes} $\{[\CO_{L}^{w}] \mid w \in W^{L}\}$ is an $R(T)$-basis of $K_{T}(X)$, and hence it is an $R(T)\bra{Q_{i_{1}}, \ldots, Q_{i_{k}}}$-basis of $QK_{T}(X)$. If $L = \emptyset$, then we write as $\CO^{w}_{L} = \CO^{w}$. 

\subsection{\texorpdfstring{Hirzebruch's $\lambda_{y}$-classes}{Hirzebruch's lambda-classes}}
Let $\CE$ be a $T$-equivariant vector bundle over $X$ of rank $r$. 
Then the \emph{Hirzebruch's $\lambda_{y}$-class} $\lambda_{y}(\CE) \in K_{T}(X)[y]$, introduced in \cite{Hirzebruch}, is defined by 
\begin{equation}
\lambda_{y}(\CE) := \sum_{d=0}^{r} y^{d} \class{\extprod^{d} \CE}. 
\end{equation}
For an exact sequence $0 \rightarrow \CE_{1} \rightarrow \CE_{2} \rightarrow \CE_{3} \rightarrow 0$, there is a well-known relation 
\begin{equation}
\lambda_{y}(\CE_{2}) = \lambda_{y}(\CE_{1}) \cdot \lambda_{y}(\CE_{3})
\end{equation} 
in $K_{T}(X)[y]$, where $\cdot$ denotes the (classical) product of $K_{T}(X)$. This relation is called a \emph{Whitney relation}. 

\subsection{\texorpdfstring{Type $C$ setting}{Type C setting}}
Consider the case where $G$ is of type $C_{n}$. Note that we use the same setting as \cite[Section~2.4]{Kouno}.
Let $J_{n}$ be the $n \times n$-matrix defined by 
\begin{equation}
J_{n} := \begin{bmatrix}
0 & 0 & \cdots & 0 & 1 \\ 
0 & 0 & \cdots & 1 & 0 \\ 
\vdots & \vdots &  & \vdots & \vdots \\ 
0 & 1 & \cdots & 0 & 0 \\ 
1 & 0 & \cdots & 0 & 0
\end{bmatrix}. 
\end{equation}
Set 
\begin{equation}
J_{2n}^{-} := \begin{bmatrix}
O & J_{n} \\ 
-J_{n} & O 
\end{bmatrix}. 
\end{equation}
Then the non-degenerate skew-symmetric form $(-, -)$ on $\BC^{2n}$ is defined as $(v, w) := {}^{t}v J_{2n}^{-} w$ for $v, w \in \BC^{2n}$. 
The \emph{symplectic group} $\Sp{2n}$ is defined by 
\begin{equation}
\Sp{2n} := \{ A \in M_{2n}(\BC) \mid {}^{t}A J_{2n}^{-} A = J_{2n}^{-} \}; 
\end{equation}
here $M_{2n}(\BC)$ denotes the ring of all $2n \times 2n$-matrices with complex entries. 
Then $\Sp{2n}$ is an algebraic group of type $C_{n}$ and its Dynkin diagram is displayed as follows: 
\begin{equation}
\dynkin[labels={1,2,n-1,n}, root radius=0.1cm, edge length=1.5cm, text style/.style={scale=1}]{C}{oo.oo}
\end{equation}
We take a Borel subgroup $B$ of $\Sp{2n}$ as the subgroup of $\Sp{2n}$ composed of all upper triangular matrices. 
The \emph{flag manifold} $\Sp{2n}/B$ of type $C_{n}$ is isomorphic to the following manifold
\begin{equation}
\FlC{n} := \left\{ (0 \subset F_{1} \subset \cdots \subset F_{n} \subset \BC^{2n}) \ \middle| \ \parbox{22em}{for $1 \le i \le n$, $F_{i}$ is an isotropic subspace of $\BC^{2n}$ with respect to $(-, -)$ with $\dim_{\BC}F_{i} = i$} \right\}. 
\end{equation}
We can take a tautological filtration $0 = \CS_{0} \subset \CS_{1} \subset \cdots \subset \CS_{n} \subset \udBC^{2n}$ of the trivial bundle $\udBC^{2n}$. 

Also, we consider the partial flag manifold $\Sp{2n}/P_{L}$ in the case $L = I \setminus \{n\}$. In this case, this is isomorphic to the \emph{Lagrangian Grassmannian} defined as 
\begin{equation}
\LG(n) := \{ F \mid \text{$F$ is an isotropic subspace of $\BC^{2n}$ with respect to $(-, -)$ with $\dim_{\BC}F = n$} \}. 
\end{equation}
We can take the tautological bundle $\CS$ and the universal quotient bundle $\CQ$ on $\LG(n)$. Note that $\CQ \simeq \CS^{\vee}$. 

Our maximal torus $T$ is written as 
\begin{equation}
T = \{ \diag(t_{1}, t_{2}, \ldots, t_{n}, t_{n}^{-1}, \ldots, t_{2}^{-1}, t_{1}^{-1}) \mid t_{1}, \ldots, t_{n} \in \BC^{\ast} = \BC \setminus \{0\} \}; 
\end{equation}
here $\diag(x_{1}, \ldots, x_{2n})$ denotes the $2n \times 2n$-diagonal matrix whose entries are $x_{1}, \ldots, x_{2n}$. 
Its Lie algebra $\Fh = \Lie(T)$ is written as 
\begin{equation}
\Fh = \{ \diag (a_{1}, a_{2}, \ldots, a_{n}, -a_{n}, \ldots, -a_{2}, -a_{1}) \mid a_{1}, \ldots, a_{n} \in \BC \}. 
\end{equation}
For $1 \le i \le n$, define $\ve_{i}: \Fh \rightarrow \BC$ by 
\begin{equation}
\ve_{i}(\diag(a_{1}, a_{2}, \ldots, a_{n}, -a_{n}, \ldots, -a_{2}, -a_{1})) := a_{i}. 
\end{equation}
We see that $\Lambda = \bigoplus_{1 \le i \le n} \BZ\ve_{i}$ and $R(T) = \BZ[\e^{\pm \ve_{1}}, \ldots, \e^{\pm \ve_{n}}]$. 
Fundamental weights $\vpi_{i}$, $1 \le i \le n$, are written as $\vpi_{i} = \ve_{1} + \cdots + \ve_{i}$. 

For simplicity, set $[a, b] := \{a < a+1 < \cdots < b\}$ for $1 \le a < b \le n$. 
For $J \subset [1, n]$, define 
\begin{equation}
\ve_{J} := \sum_{j \in J} \ve_{j}. 
\end{equation}
Then we can write as 
\begin{equation}
\class{\extprod^{d} \CS_{k}} = \sum_{\substack{J \subset [1, k] \\ |J| = d}} [\CO(-\ve_{J})], \quad \class{\extprod^{d} \CS_{k}^{\vee}} = \sum_{\substack{J \subset [1, k] \\ |J| = d}} [\CO(\ve_{J})]. 
\end{equation}

\subsection{\texorpdfstring{$QK$-Whitney relations for type $C$ flag manifolds}{QK-Whitney relations for type C}}
The quantum $K$-theoretic Whitney relations in $QK_{T}(\Fl^{C_{n}})$, described in \cite{Kouno}, is written as follows: 
\begin{thm}[{\cite[Theorem~3.2]{Kouno}}] \label{thm:QK-Whitney_full}
In $QK_{T}(\Fl^{C_{n}})$, we have 
\begin{equation} \label{eq:QK-Whitney_full}
\begin{split}
& \lambda_{y}(\CS_{n}) \star \lambda_{y}(\CS_{n}^{\vee}) = \lambda_{y}(\udBC^{2n}) \\ 
& \hspace{10mm} - y^{2} \sum_{p=1}^{n} \frac{Q_{p}Q_{p+1} \cdots Q_{n}}{1-Q_{p-1}} (\lambda_{y}(\CS_{p-1}) - Q_{p-1} \lambda_{y}(\CS_{p-2})) \star (\lambda_{y}(\CS_{p-1}^{\vee}) - Q_{p-1} \lambda_{y}(\CS_{p-2}^{\vee})); 
\end{split}
\end{equation}
here, we understand that $\CS_{-1} = \CS_{-1}^{\vee} = 0$ and $Q_{0} = 0$. 
\end{thm}

\begin{rem}
Theorem~\ref{thm:QK-Whitney_full} does not give the complete list of $QK$-Whitney relations; in fact, there are other relations. 
Since we do not use those relations in this paper, we omit the review of those. For details, see \cite[Theorem~3.1]{Kouno}. 
\end{rem}

In \cite{Kouno}, the author proved \eqref{eq:QK-Whitney_full} by using the following lemma, which we will use again in this paper. 

\begin{lem}[{\cite[Lemma~5.7]{Kouno}}] \label{lem:to_recursive}
In $QK_{T}(\FlC{n})$, for $1 \le d \le 2n$ and $1 \le p \le n$, we have
\begin{equation} \label{eq:QK-Whitney_computation_lem2}
\begin{split}
& \sum_{\substack{J, K \subset [1, n] \\ |J| + |K| = d \\ \max J = \max K = p}} [\CO(-\ve_{J})] \star [\CO(\ve_{K})] \\ 
& \qquad = \frac{1-Q_{p}}{1-Q_{p-1}} \sum_{\substack{k, l \ge 0 \\ k+l = d-2}} \left( \left[ \extprod^{k} \CS_{p-1} \right] - Q_{p-1} \left[ \extprod^{k} \CS_{p-2} \right] \right) \\ 
& \hspace{60mm} \star \left( \left[ \extprod^{l} \CS_{p-1}^{\vee} \right] - Q_{p-1} \left[ \extprod^{l} \CS_{p-2}^{\vee} \right] \right). 
\end{split}
\end{equation}
\end{lem}

\section{Kato's pushforward and its explicit computation} \label{sec:Kato_push}
First, we recall Kato's pushforward $QK_{T}(\FlC{n}) \twoheadrightarrow QK_{T}(\LG(n))$. 
Then, to compute images under Kato's pushforward, we review quantum alcove models and Chevalley formulas. 
After that, we actually compute images of some line bundles under Kato's pushforward. 

\subsection{Kato's pushforward}
Let $G$ be of arbitrary type in this section. Take $L \subset I$. 
There exists a canonical projection $\pi: G/B \twoheadrightarrow G/P_{L}$. Kato established the pushforward $\pi_{\ast}: QK_{T}(G/B) \twoheadrightarrow QK_{T}(G/P_{L})$. 
\begin{thm}[{\cite[Theorem~2.19]{Kato}}]
There exists an $R(T)$-algebra surjection $\pi_{\ast}: QK_{T}(G/B) \twoheadrightarrow QK_{T}(G/P_{L})$ such that 
$\pi_{\ast}(Q^{\beta} \CO^{w}) = Q^{[\beta]} \CO_{L}^{\mcr{w}}$ for $\beta \in Q^{\vee, +}$ and $w \in W$; here $\mcr{w} \in W^{L}$ denotes the minimal-length representative of the coset $wW_{L} \in W/W_{L}$ for $w \in W$ and $[\beta] := \sum_{j \in I \setminus L} c_{j}\alpha_{j}^{\vee} \in Q_{L}^{\vee, +}$ for $\beta = \sum_{j \in I} c_{j}\alpha_{j}^{\vee} \in Q^{\vee, +}$. 
\end{thm}

Note that $\pi_{\ast}$ can be extended to the surjection $QK_{T}(G/B)[y] \twoheadrightarrow QK_{T}(G/P_{L})[y]$ canonically, for which we use the same notation $\pi_{\ast}$. 
Also, the restriction $\pi_{\ast}|_{K_{T}(G/B)}$ coincides with the classical pushforward $\pi_{!}: K_{T}(G/B) \rightarrow K_{T}(G/P)$; note that $\pi_{!}$ is not a ring homomorphism. 

In what follows, we use this pushforward in the case $G = \Sp{2n}$ and $L = I \setminus \{n\}$; then $G/B = \FlC{n}$ and $G/P_{L} = \LG(n)$, i.e., we consider the pushforward $\pi_{\ast}: QK_{T}(\FlC{n}) \twoheadrightarrow QK_{T}(\LG(n))$. 
In this situation, we write $Q = Q_{n}$ in $QK_{T}(\LG(n))$. 
For $1 \le d \le n$, we have 
\begin{equation}
\pi_{\ast}\left( \class{\extprod^{d} \CS_{n}} \right) = \class{\extprod^{d} \CS}, \quad \pi_{\ast}\left( \class{\extprod^{d} \CS_{n}^{\vee}} \right) = \class{\extprod^{d} \CQ} 
\end{equation}
since $\class{\extprod^{d} \CS_{n}}$ (resp., $\class{\extprod^{d} \CS_{n}^{\vee}}$) coincides with $\pi^{\ast} \left(\class{\extprod^{d} \CS}\right)$ (resp., $\pi^{\ast} \left(\class{\extprod^{d} \CQ} \right)$), where $\pi^{\ast}: K_{T}(\LG(n)) \rightarrow K_{T}(\FlC{n})$ denotes the pullback of $\pi$. 
To compute the image of the right-hand side of \eqref{eq:QK-Whitney_full} under $\pi_{\ast}$, we show the following lemma in the next section. 
\begin{lem} \label{lem:line_proj}
\begin{enu}
\item We have $\pi_{\ast}([\CO(\ve_{J})]) = 0$ if $J \not= [1, k]$ for any $1 \le k \le n$. 

\item We have $\pi_{\ast}([\CO(-\ve_{J})]) = 0$ if $J \not= [k, n]$ for any $1 \le k \le n$.
\end{enu}
\end{lem}

As a corollary, we obtain the following two identities. 
\begin{cor} \label{cor:quotient}
For $1 \le d \le n$, we have 
\begin{equation}
\class{\extprod^{d} \CQ} = \pi_{\ast}(\class{\CO(\ve_{1} + \cdots + \ve_{d})}) 
\end{equation}
in $QK_{T}(\LG(n))$. 
\end{cor}
\begin{proof}
By Lemma~\ref{lem:line_proj}\,(1), we obtain 
\begin{equation}
\class{\extprod^{d} \CQ} = \pi_{\ast} \left( \class{\extprod^{d} \CS_{n}^{\vee}} \right) = \sum_{\substack{J \subset [1, n] \\ |J|=d}} \pi_{\ast}(\class{\CO(\ve_{J})}) = \pi_{\ast}\left( \class{\CO(\ve_{[1, d]})} \right),  
\end{equation}
as desired. 
\end{proof}

\begin{cor} \label{cor:proj_lambda}
For $1 \le k \le n-1$, we have $\pi_{\ast}(\lambda_{y}(\CS_{k})) = 1$. 
\end{cor}
\begin{proof}
Since $J$ can not contain $n$ if $J \subset [1, k]$, Lemma~\ref{lem:line_proj}\,(2) implies
\begin{equation}
\pi_{\ast}(\lambda_{y}(\CS_{k})) = \sum_{d=0}^{k} y^{d} \sum_{\substack{J \subset [1, k] \\ |J|=d}} \pi_{\ast}(\class{\CO(-\ve_{J})}) = 1, 
\end{equation}
as desired. 
\end{proof}

\subsection{Quantum alcove model and Chevalley formula}
To use Kato's pushforward, we need an expansion of line bundles into linear combinations of Schubert classes. 
For this purpose, we use the $K_{T}$-Chevalley formula described explicitly by Lenart and Postnikov \cite{LP}. 
Since $K_{T}$-Chevalley formulas are written in terms of alcove models, we first recall their definition. 
We also use \emph{quantum} alcove models introduced by Lenart and Lubovsky \cite{LL} (and developed by Lenart, Naito, and Sagaki \cite{LNS}). For a while, we consider the case that $G$ is of arbitrary type. 

\begin{defn}[{\cite[Definition~6.1]{BFP}}]
The \emph{Quantum Bruhat graph} $\QBG(W)$ is a labeled directed graph whose vertex set is $W$, and for $x, y \in W$ and $\alpha \in \Delta^{+}$, we have an edge $x \xrightarrow{\alpha} y$ if $y = xs_{\alpha}$ and one of the following holds: 
\begin{itemize}
\item[(B)] $\ell(y) = \ell(x) + 1$, or 
\item[(Q)] $\ell(y) = \ell(x) - 2\pair{\rho}{\alpha^{\vee}} + 1$. 
\end{itemize}
If (B) (resp., (Q)) holds, then the corresponding edge is called a \emph{Bruhat edge} (resp., \emph{quantum edge}). 
\end{defn}

Set $\Fh_{\BR}^{\ast} := \Lambda \otimes_{\BZ} \BR$. For $\alpha \in \Delta$ and $k \in \BZ$, define a hyperplane $H_{\alpha, k}$ of $\Fh_{\BR}^{\ast}$ by 
\begin{equation}
H_{\alpha, k} := \{ \xi \in \Fh_{\BR}^{\ast} \mid \pair{\xi}{\alpha^{\vee}} = k \}. 
\end{equation}
An \emph{alcove} is a connected component of $\Fh_{\BR}^{\ast} \setminus \bigcup_{\alpha \in \Delta, k \in \BZ} H_{\alpha, k}$. Then an alcove $A_{\circ}$ defined as 
\begin{equation}
A_{\circ} := \{ \xi \in \Fh_{\BR}^{\ast} \mid \text{$0 < \pair{\xi}{\alpha^{\vee}} < 1$ for $\alpha \in \Delta^{+}$} \}
\end{equation}
is called the \emph{fundamental alcove}. 
For $\lambda \in \Lambda$, define 
\begin{equation}
A_{\lambda} := A_{\circ} + \lambda = \{ \xi + \lambda \mid \xi \in A_{\circ} \}. 
\end{equation} 
Alcoves $A$ and $B$ are called \emph{adjacent} if $A$ and $B$ share a common wall. 
For alcoves $A$, $B$, and $\alpha \in \Delta$, we write $A \xrightarrow{\alpha} B$ if
\begin{itemize}
\item $A$ and $B$ are adjacent, 
\item the common wall of $A$ and $B$ is contained in $H_{\alpha, k}$ for some $k \in \BZ$, and 
\item $\alpha$ points the direction from $A$ to $B$. 
\end{itemize}

\begin{defn}[{\cite[Definitions~5.2 and 5.4]{LP}}]
\begin{enu}
\item A sequence $(A_{0}, A_{1}, \ldots, A_{r})$ is called an \emph{alcove path} if $A_{k-1}$ and $A_{k}$ are adjacent for $1 \le k \le r$. 
\item An alcove path $(A = A_{0}, \ldots, A_{r} = B)$ is called \emph{reduced} if $r$ is minimal among all alcove paths from $A$ to $B$. 
\item Let $\lambda \in \Lambda$. A sequence $\Gamma = (\gamma_{1}, \ldots, \gamma_{r})$ is called a \emph{$\lambda$-chain} if there exists an alcove path $A^{\Gamma} = (A_{\circ} = A_{0}, A_{1}, \ldots, A_{r} = A_{-\lambda})$ such that $A_{0} \xrightarrow{-\gamma_{1}} A_{1} \xrightarrow{-\gamma_{2}} \cdots \xrightarrow{-\gamma_{r}} A_{r}$. In addition, we say $\Gamma$ is \emph{reduced} if the corresponding alcove path $A^{\Gamma}$ is reduced. 
\end{enu}
\end{defn}

We introduce objects of the theory of quantum alcove models, called \emph{admissible subsets}. 

\begin{defn}[{\cite[Definition~3.4]{LL} and \cite[Definition~17]{LNS}}]
Let $\lambda \in \Lambda$ and take a $\lambda$-chain $\Gamma = (\gamma_{1}, \ldots, \gamma_{r})$. Also, let $w \in W$. 
\begin{enu}
\item A set $A = \{ i_{1} < i_{2} < \cdots < i_{t} \} \subset \{1, \ldots, r\}$ is called \emph{$w$-admissible} if there exists a directed path 
\begin{equation}
\Pi(w, A): w = w_{0} \xrightarrow{|\gamma_{i_{1}}|} w_{1} \xrightarrow{|\gamma_{i_{2}}|} \cdots \xrightarrow{|\gamma_{i_{t}}|} w_{t} 
\end{equation}
in $\QBG(W)$. 
We denote by $\CA(w, \Gamma)$ the set of all $w$-admissible subsets associated to $\Gamma$. 

\item Set 
\begin{align}
\ed(A) &:= w_{t}, \\ 
\down(A) &:= \sum_{\substack{1 \le k \le t \\ \text{$w_{k-1} \rightarrow w_{k}$ is a quantum edge}}} |\gamma_{i_{k}}|^{\vee}, \\ 
n(A) &:= \# \{ j \in A \mid \gamma_{j} \in -\Delta^{+} \}. 
\end{align}
\item Let $(A_{\circ} = A_{0}, A_{1}, \ldots, A_{r} = A_{-\lambda})$ be the alcove path corresponding to $\Gamma$. For $1 \le k \le r$, we denote by $\wh{r}_{k}: \Fh_{\BR}^{\ast} \rightarrow \Fh_{\BR}^{\ast}$ the reflection with respect to the hyperplane containing the common wall of $A_{k-1}$ and $A_{k}$. For $A = \{i_{1} < \cdots < i_{t}\} \in \CA(w, \Gamma)$, set 
\begin{equation}
\wt(A) := -w\wh{r}_{i_{1}} \cdots \wh{r}_{i_{t}} (-\lambda). 
\end{equation}

\item For $w \in W$ and a $\lambda$-chain $\Gamma$ associated to $\lambda \in \Lambda$, set 
\begin{equation}
\CA|_{Q=0}(w, \Gamma) := \{ A \in \CA(w, \Gamma) \mid \down(A) = 0 \}. 
\end{equation}
\end{enu}
\end{defn}

The classical \emph{$K_{T}$-Chevalley formula} is written explicitly as follows. 
\begin{thm}[{\cite[Theorem~6.1]{LP}; see also \cite[Remark~4.9]{Kouno}}] \label{thm:Chevalley}
Let $\lambda \in \Lambda$ and take an arbitrary reduced $\lambda$-chain $\Gamma$. 
Let $w \in W$. 
Then in $K_{T}(G/B)$, we have 
\begin{equation}
[\CO(\lambda)] \cdot [\CO^{w}] = \sum_{A \in \CA|_{Q=0}(w, \Gamma)} (-1)^{n(A)} e^{{}-\wt(A)} [\CO^{\ed(A)}]. 
\end{equation}
\end{thm}
\subsection{Proof of Lemma~\ref{lem:line_proj}}
We give a proof of Lemma~\ref{lem:line_proj}. In this section, we consider the case $G = \Sp{2n}$ and $L = I \setminus \{n\}$ again. 
Recall an explicit construction of a reduced $\ve_{J}$-chain $\Gamma_{\vtl}(\ve_{J})$ and a reduced $(-\ve_{J})$-chain $\Gamma(-\ve_{J})$ from \cite[Section~A.2]{Kouno}. 
In the case that $J \not= [1, k]$ for any $1 \le k \le n$, there exist unique $1 \le d' \le n-2$ and $1 \le d \le n-2$ such that $d', d'+1, \ldots, d \in J$, $1, \ldots, d' \notin J$ and $d+1 \notin J$. Note that there exists some $p > d+1$ such that $p \in J$. In this case, if we write $\Gamma_{\vtl}(\ve_{J}) = (\gamma_{1}, \ldots, \gamma_{l})$, then $\gamma_{l} = -\alpha_{d+1}$. 
Since there always exists an edge $w \xrightarrow{\alpha} ws_{\alpha}$ in $\QBG(W)$ for all $w \in W$ and all simple root $\alpha$, we can construct a sign-reversing involution $\iota: \CA(e, \Gamma_{\vtl}(\ve_{J})) \xrightarrow{\sim} \CA(e, \Gamma_{\vtl}(\ve_{J}))$ defined as
\begin{equation}
\iota(A) := \begin{cases} A \setminus \{l\} & \text{if $l \in A$,} \\ A \sqcup \{l\} & \text{if $l \notin A$.} \end{cases}
\end{equation}
By the definition of $\iota$, we have $\ed(\iota(A)) = \ed(A) s_{d+1}$ for $A \in \CA(e, \Gamma_{\vtl}(\ve_{J}))$. 
Also, we can see that the common wall of $A_{l-1}$ and $A_{l}$ is contained in $H_{\alpha_{d+1}, 1}$. Since $\pair{\ve_{J}}{\alpha_{d+1}} = -1$, we have $\wh{r}_{l}(-\ve_{J}) = -\ve_{J}$. Hence we have $\wt(\iota(A)) = \wt(A)$ for $A \in \CA(e, \Gamma_{\vtl}(\ve_{J}))$. 

In addition, we have the following ``down-free'' lemma. 
\begin{lem}[{\cite[Proposition~A.3]{Kouno}}] \label{lem:down-free}
\begin{enu}
\item For $A \in \CA(e, \Gamma_{\vtl}(\ve_{J}))$, we have $\down(A) = 0$. 
\item For $A \in \CA(e, \Gamma(-\ve_{J}))$, we have $\down(A) = 0$. 
\end{enu}
\end{lem}

\begin{proof}[Proof of Lemma~\ref{lem:line_proj}]
We only give a proof of (1); (2) can be shown by the same argument. 
Assume, as above, that there exist unique $1 \le d' \le n-2$ and $1 \le d \le n-2$ such that $d', d'+1, \ldots, d \in J$, $1, \ldots, d'-1 \notin J$ and $d+1 \notin J$. 
By applying Theorem~\ref{thm:Chevalley} for $\lambda = \ve_{J}$ and $w = e$, combining with Lemma~\ref{lem:down-free}, we obtain 
\begin{equation}
[\CO(\ve_{J})] = \sum_{A \in \CA(e, \Gamma_{\vtl}(\ve_{J}))} (-1)^{n(A)} e^{{}-\wt(A)} \CO^{\ed(A)}. 
\end{equation}
Set 
\begin{equation}
\CA_{0}(e, \Gamma_{\vtl}(\ve_{J})) := \{ A \in \CA(e, \Gamma_{\vtl}(\ve_{J})) \mid l \notin A \}. \\ 
\end{equation}
Then by using an involution $\iota$, we obtain 
\begin{equation}
[\CO(\ve_{J})] = \sum_{A \in \CA_{0}(e, \Gamma_{\vtl}(\ve_{J}))} (-1)^{n(A)} e^{{}-\wt(A)} (\CO^{\ed(A)} - \CO^{\ed(A)s_{d+1}}). 
\end{equation}
Since $\mcr{ws_{d+1}} = \mcr{w}$ for $w \in W$, we have $\pi_{\ast}([\CO(\ve_{J})]) = 0$, as desired. 
\end{proof}

\section{Recursive identities} \label{sec:recursive}
We derive key recursive identities in $QK_{T}(\LG(n))[y]$, and then we give a proof of Theorem~\ref{thm:Whitney_rel}. 
Also, we give examples of the computation of $QK$-Whitney relations for type $C_{2}$ and $C_{3}$. 

\subsection{\texorpdfstring{Recursive identities and proof of $QK$-Whitney relations}{Recursive identities and proof of QK-Whitney relations}}
We describe a recursive identity in $QK_{T}(\LG(n))[y]$, which plays a crucial role in the proof of Theorem~\ref{thm:Whitney_rel}. 
For $0 \le p \le n$, define $\CP_{p} \in QK_{T}(\FlC{n})[y]$ by 
\begin{equation} \label{eq:def_Pp}
\CP_{p} := \frac{1}{1-Q_{p}}(\lambda_{y}(\CS_{p}) - Q_{p}\lambda_{y}(\CS_{p-1})) \star (\lambda_{y}(\CS_{p}^{\vee}) - Q_{p}\lambda_{y}(\CS_{p-1}^{\vee}));  
\end{equation}
here, we understand that $\CS_{-1} = \CS_{-1}^{\vee} = 0$ and $Q_{0} = 0$. 

\begin{lem}
For $1 \le p \le n-1$, we have 
\begin{equation} \label{eq:recursive}
\pi_{\ast}(\CP_{p}) = y^{2} \pi_{\ast}(\CP_{p-1}) + y^{p} \class{\extprod^{p} \CQ}. 
\end{equation}
\end{lem}
\begin{proof}
Observe that 
\begin{equation}
\lambda_{y}(\CS_{p}) = \lambda_{y}(\CS_{p-1}) + \sum_{d=0}^{p-1} y^{d+1} \sum_{\substack{J \subset [1, p-1] \\ |J| = d}} \class{\CO(-\ve_{J \sqcup \{p\}})}
\end{equation}
by the definition of $\lambda_{y}$-classes. Similarly, we see that 
\begin{equation}
\lambda_{y}(\CS_{p}^{\vee}) = \lambda_{y}(\CS_{p-1}^{\vee}) + \sum_{d=0}^{p-1} y^{d+1} \sum_{\substack{J \subset [1, p-1] \\ |J| = d}} \class{\CO(\ve_{J \sqcup \{p\}})}. 
\end{equation}
Hence we compute as 
\begin{align}
\begin{split}
\CP_{p} &= \frac{1}{1-Q_{p}} \left( (1-Q_{p}) \lambda_{y}(\CS_{p-1}) + \sum_{d=0}^{p-1} y^{d+1} \sum_{\substack{J \subset [1, p-1] \\ |J| = d}} \class{\CO(-\ve_{J \sqcup \{p\}})} \right) \\ 
& \hspace{30mm} \star \left( (1-Q_{p}) \lambda_{y}(\CS_{p-1}^{\vee}) + \sum_{d=0}^{p-1} y^{d+1} \sum_{\substack{J \subset [1, p-1] \\ |J| = d}} \class{\CO(\ve_{J \sqcup \{p\}})} \right) 
\end{split} \\ 
\begin{split}
&= \underbrace{(1-Q_{p}) \lambda_{y}(\CS_{p-1}) \star \lambda_{y}(\CS_{p-1}^{\vee})}_{:= \CP_{p}^{(1)}} \\ 
& \qquad + \underbrace{\lambda_{y}(\CS_{p-1}) \star \sum_{d=0}^{p-1} y^{d+1} \sum_{\substack{J \subset [1, p-1] \\ |J| = d}} \class{\CO(\ve_{J \sqcup \{p\}})}}_{:= \CP_{p}^{(2)}} \\ 
& \qquad + \underbrace{\lambda_{y}(\CS_{p-1}^{\vee}) \star \sum_{d=0}^{p-1} y^{d+1} \sum_{\substack{J \subset [1, p-1] \\ |J| = d}} \class{\CO(-\ve_{J \sqcup \{p\}})}}_{:= \CP_{p}^{(3)}} \\ 
& \qquad + \underbrace{\frac{1}{1-Q_{p}}\sum_{d=0}^{p-1} \sum_{d'=0}^{p-1} y^{d+d'+2} \sum_{\substack{J \subset [1, p-1] \\ |J|=d}} \sum_{\substack{K \subset [1, p-1] \\ |K|=d'}} \class{\CO(-\ve_{J\sqcup\{p\}})} \star \class{\CO(\ve_{J\sqcup\{p\}})}}_{:= \CP_{p}^{(4)}}. 
\end{split}
\end{align}
Since $\pi_{\ast}(Q_{p}) = 1$, we have $\pi_{\ast}(\CP_{p}^{(1)}) = 0$. Also, by Lemma~\ref{lem:line_proj}\,(1) and Corollaries~\ref{cor:quotient} and \ref{cor:proj_lambda}, we have 
\begin{equation}
\pi_{\ast}(\CP_{p}^{(2)}) = 1 \star \sum_{d=0}^{p-1}y^{d+1} \sum_{\substack{J \subset [1, p-1] \\ |J|=d}} \pi_{\ast}\left( \class{\CO(\ve_{J\sqcup\{p\}})} \right) = y^{p} \pi_{\ast}\left( \class{\CO(\ve_{[1,p]})} \right) = y^{p} \class{\extprod^{p} \CQ}. 
\end{equation}
In addition, we see that 
\begin{equation}
\pi_{\ast}(\CP_{p}^{(3)}) = \pi_{\ast}(\lambda_{y}(\CS_{p-1}^{\vee})) \star \sum_{d=0}^{p-1} y^{d+1} \sum_{\substack{J \subset [1, p-1] \\ |J|=d}} \pi_{\ast} \left( \class{\CO(-\ve_{J\sqcup\{p\}})} \right) = 0
\end{equation}
since $J \sqcup \{p\}$ can not contain $n$ for all $J \subset [1, p-1]$. To consider $\pi_{\ast}(\CP_{p}^{(4)})$, observe that 
\begin{align}
\CP_{p}^{(4)} &= \frac{1}{1-Q_{p}}y^{2} \sum_{j=0}^{2(p-1)} y^{j} \sum_{\substack{J, K \subset [1, p-1] \\ |J|+|K| = j}} \class{\CO(-\ve_{J\sqcup\{p\}})} \star \class{\CO(-\ve_{K\sqcup\{p\}})} \\ 
&= \frac{1}{1-Q_{p}} y^{2} \sum_{j=0}^{2(p-1)} y^{j} \sum_{\substack{J, K \subset [1, n] \\ |J|+|K| = j+2 \\ \max J = \max K = p}} \class{\CO(-\ve_{J})} \star \class{\CO(\ve_{K})} \\ 
\begin{split}
&= \frac{1}{1-Q_{p}} y^{2} \sum_{j=0}^{2(p-1)} y^{j} \frac{1-Q_{p}}{1-Q_{p-1}} \sum_{\substack{k, l \ge 0 \\ k+l = j}} \left( \left[ \extprod^{k} \CS_{p-1} \right] - Q_{p-1} \left[ \extprod^{k} \CS_{p-2} \right] \right) \\ 
& \hspace{80mm} \star \left( \left[ \extprod^{l} \CS_{p-1}^{\vee} \right] - Q_{p-1} \left[ \extprod^{l} \CS_{p-2}^{\vee} \right] \right) 
\end{split} \\ 
&= \frac{1}{1-Q_{p}} \frac{1-Q_{p}}{1-Q_{p-1}} y^{2} (\lambda_{y}(\CS_{p-1}) - Q_{p-1}\lambda_{y}(\CS_{p-2})) \star (\lambda_{y}(\CS_{p-1}^{\vee}) - Q_{p-1} \lambda_{y}(\CS_{p-2}^{\vee})) \\ 
&= y^{2} \CP_{p-1}; 
\end{align} 
here, in the third equality, we used Lemma~\ref{lem:to_recursive}.  
Therefore, we obtain 
\begin{equation}
\pi_{\ast}(\CP_{p}) = \pi_{\ast}(\CP_{p}^{(1)}) + \pi_{\ast}(\CP_{p}^{(2)}) + \pi_{\ast}(\CP_{p}^{(3)}) + \pi_{\ast}(\CP_{p}^{(4)}) = y^{p} \class{\extprod^{p} \CQ} + y^{2}\pi_{\ast}(\CP_{p-1}), 
\end{equation}
as desired. 
\end{proof}

\begin{proof}[Proof of Theorem~\ref{thm:Whitney_rel}]
By applying $\pi_{\ast}: QK_{T}(\FlC{n})[y] \twoheadrightarrow QK_{T}(\LG(n))[y]$ to \eqref{eq:QK-Whitney_full}, we obtain 
\begin{align}
\lambda_{y}(\CS) \star \lambda_{y}(\CQ) &= \lambda_{y}(\udBC^{2n}) - y^{2} \sum_{p=1}^{n} \pi_{\ast}(Q_{p}Q_{p+1} \cdots Q_{n}) \pi_{\ast}(\CP_{p-1}) \\ 
&= \lambda_{y}(\udBC^{2n}) - y^{2} Q \sum_{p=1}^{n} \pi_{\ast}(\CP_{p-1}). 
\end{align} 
Since $\CP_{0} = 1$, \eqref{eq:recursive} implies
\begin{equation} \label{eq:pi_ast_P}
\pi_{\ast}(\CP_{p-1}) = y^{p-1} \class{\extprod^{p-1} \CQ} + y^{p} \class{\extprod^{p-2} \CQ} + \cdots + y^{2p-2}. 
\end{equation}
Hence, in the sum $\sum_{p=1}^{n} \pi_{\ast}(\CP_{p-1})$, the coefficient of $y^{d}$ for $0 \le d \le 2n-2$ is 
\begin{equation}
\sum_{p = \lceil d/2 \rceil + 1}^{\min\{d+1, n\}} \class{\extprod^{2p-d-2} \CQ}, 
\end{equation}
which coincides with $\CF_{d}$. 
This proves the theorem. 
\end{proof}

\subsection{\texorpdfstring{Examples of $QK$-Whitney relations}{Examples of QK-Whitney relations}}
We confirm Theorem~\ref{thm:Whitney_rel} in type $C_{2}$ by computing the image of the right-hand side of \eqref{eq:QK-Whitney_full} under $\pi_{\ast}$ directly. 
\begin{exm}
Let $n = 2$ and consider $QK_{T}(\LG(2))$. 
Then 
\begin{align}
\CP_{0} &= 1, \\ 
\CP_{1} &= \frac{1}{1-Q_{1}} (\lambda_{y}(\CS_{1}) - Q_{1}) \star (\lambda_{y}(\CS_{1}^{\vee}) - Q_{1}) \\ 
&= \frac{1}{1-Q_{1}} (1 + y\class{\CO(-\ve_{1})} - Q_{1}) \star (1 + y\class{\CO(\ve_{1})} - Q_{1}) \\ 
&= \frac{1}{1-Q_{1}} ((1-Q_{1})^{2} + y(1-Q_{1})(\class{\CO(-\ve_{1})} + \class{\CO(\ve_{1})}) + y^{2} \underbrace{\class{\CO(-\ve_{1})} \star \class{\CO(\ve_{1})}}_{\mathrlap{\hspace{-14mm}= 1-Q_{1} \text{ by \cite[Lemma~6.7]{KN}}}}) \\
&= (1-Q_{1}) + y(\class{\CO(-\ve_{1})} + \class{\CO(\ve_{1})}) + y^{2}. 
\end{align}
Hence 
\begin{align}
\lambda_{y}(\CS) \star \lambda_{y}(\CQ) &= \lambda_{y}(\udBC^{4}) - y^{2} (\pi_{\ast}(Q_{1}Q_{2}\CP_{0}) + \pi_{\ast} (Q_{2}\CP_{1})) \\ 
&= \lambda_{y}(\udBC^{4}) - y^{2} Q (1 + y (\underbrace{\pi_{\ast}(\class{\CO(-\ve_{1})})}_{\mathrlap{\hspace{-10mm}=0 \text{ by Lemma~\ref{lem:line_proj}}}} + \underbrace{\pi_{\ast}(\class{\CO(\ve_{1})})}_{\mathrlap{\hspace{-8mm}= \class{\CQ} \text{ by Corollary~\ref{cor:quotient}}}}) + y^{2}) \\ 
&= \lambda_{y}(\udBC^{4}) - y^{2}Q (1 + y \class{\CQ} + y^{2}). 
\end{align}
This coincides with \eqref{eq:QK-Whitney_rel} because $\CF_{0} = \CF_{2} = 1$, and $\CF_{1} = \class{\CQ}$.  
Note that this also coincides with \cite[(7.24)]{GMSZ}. 
\end{exm}

Next, we compute \eqref{eq:QK-Whitney_rel} for type $C_{3}$ by allowing the use of \eqref{eq:pi_ast_P}. 
\begin{exm}
Let $n = 3$ and consider $QK_{T}(\LG(3))$. By \eqref{eq:pi_ast_P}, we have 
\begin{equation}
\pi_{\ast}(\CP_{0}) = 1, \quad \pi_{\ast}(\CP_{1}) = y \class{\CQ} + y^{2}, \quad \pi_{\ast}(\CP_{2}) = y^{2}\class{\extprod^{2} \CQ} + y^{3} \class{\CQ} + y^{4}. 
\end{equation}
Therefore, we obtain 
\begin{align}
\lambda_{y}(\CS) \star \lambda_{y}(\CQ) &= \lambda_{y}(\udBC^{6}) - y^{2} Q \sum_{p=1}^{3} \pi_{\ast}(\CP_{p-1}) \\ 
&= \lambda_{y}(\udBC^{6}) - y^{2}Q \left(1 + y \class{\CQ} + y^{2}\left( 1 + \class{\extprod^{2}\CQ} \right) + y^{3} \class{\CQ} + y^{4} \right). 
\end{align}
This coincides with \cite[(7.67)]{GMSZ}. 
\end{exm}

\section{\texorpdfstring{$QK$-Whitney presentation of $QK_{T}(\LG(n))$}{QK-Whitney presentation of QKT(LG(n))}} \label{sec:presentation}
The purpose of this section is to give a proof of the $QK$-Whitney relation for $\LG(n)$. For this purpose, we review the Nakayama-type argument to obtain the $QK$-Whitney presentation from the classical one. 

\subsection{Nakayama-type argument}
To give a presentation of $QK_{T}(\LG(n))$, we recall the following Nakayama-type argument for an arbitrary $G$ and $L \subset I$, established by Gu, Mihalcea, Sharpe, Xu, Zhang, and Zou \cite{GMSXZZ2}. 

\begin{thm}[{\cite[Theorem~4.1]{GMSXZZ2}}] \label{thm:Nakayama}
Assume that we are given an isomorphism of $R(T)$-algebras
\begin{equation}
\Phi: R(T)[m_{1}, \ldots, m_{s}]/\langle P_{1}, \ldots, P_{r} \rangle \xrightarrow{\sim} K_{T}(G/P_{L})
\end{equation}
where $m_{1}, \ldots, m_{s}$ are indeterminates, and the $P_{i}$ are polynomials in the variables $m_{j}$ with coefficients in $R(T)$. 
Consider any power series
\begin{equation}
\ti{P}_{i} = \ti{P}_{i}(m_{1}, \ldots, m_{s}; Q_{1}, \ldots, Q_{k}) \in (R(T)[m_{1}, \ldots, m_{s}])\bra{Q_{1}, \ldots, Q_{k}}
\end{equation}
such that 
\begin{itemize}
\item $\ti{P}_{i}(m_{1}, \ldots, m_{s}; 0, \ldots, 0) = P_{i}(m_{1}, \ldots, m_{s})$ for all $m_{1}, \ldots, m_{s}$ and all $i$, and 
\item $\ti{P}_{i}(m_{1}, \ldots, m_{s}; Q_{1}, \ldots, Q_{k}) = 0$ in $QK_{T}(G/P_{L})$, where the $m_{i}$ are regarded as elements in $QK_{T}(G/P_{L})$ via the isomorphism $\Phi$. 
\end{itemize}
Then $\Phi$ lifts to an isomorphism of $R(T)\bra{Q_{1}, \ldots, Q_{k}}$-algebras
\begin{equation}
\ti{\Phi}: (R(T)\bra{Q_{1}, \ldots, Q_{k}})[m_{1}, \ldots, m_{s}]/\langle \ti{P}_{1}, \ldots, \ti{P}_{r} \rangle \xrightarrow{\sim} QK_{T}(G/P_{L}). 
\end{equation}
\end{thm}

\subsection{Proof of Theorem~\ref{thm:Whitney_presentation}}
We apply the Nakayama-type argument to the following Whitney presentation of $K_{T}(\LG(n))$. 
Let $R(T)[x_{1}, \ldots, x_{n}, y_{1}, \ldots, y_{n}]^{\FS_{n}}$ be the subring of the polynomial ring $R(T)[x_{1}, \ldots, x_{n}, y_{1}, \ldots, y_{n}]$ composed of all $\FS_{n}$-invariants; i.e., 
\begin{equation}
R(T)[x_{1}, \ldots, x_{n}, y_{1}, \ldots, y_{n}]^{\FS_{n}} = R(T)[e_{d}(x), e_{d}(y) \mid 1 \le d \le n].
\end{equation} 
Define $J_{\LG(n)}^{\cl}$ as an ideal of $R(T)[x_{1}, \ldots, x_{n}, y_{1}, \ldots, y_{n}]^{\FS_{n}}$ generated by 
\begin{equation}
\sum_{\substack{0 \le k, l \le d \\ k+l=d}} e_{k}(x)e_{l}(y) - e_{d}(\e^{\pm \ve}), \qquad 0 \le d \le 2n. 
\end{equation}
\begin{prop} \label{prop:classical_Whitney}
There exists an $R(T)$-algebra isomorphism
\begin{equation}
\Psi_{\cl}: R(T)[x_{1}, \ldots, x_{n}, y_{1}, \ldots, y_{n}]^{\FS_{n}}/J_{\LG(n)}^{\cl} \xrightarrow{\sim} K_{T}(\LG(n))
\end{equation}
such that $\Psi_{\cl}(e_{d}(x)) = \class{\extprod^{d} \CS}$ and $\Psi_{\cl}(e_{d}(y)) = \class{\extprod^{d} \CQ}$. 
\end{prop}
In other words, Proposition~\ref{prop:classical_Whitney} states that $K_{T}(\LG(n))$ is generated as $R(T)$-algebra by $\class{\extprod^{d} \CS}$ and $\class{\extprod^{d} \CQ}$ for $1 \le d \le n$ with defining relations deduced from the Whitney relation 
\begin{equation} \label{eq:classical_Whitney_2}
\lambda_{y}(\CS) \cdot \lambda_{y}(\CQ) = \lambda_{y}(\udBC^{2n}). 
\end{equation}
Note that the coefficient of $y^{d}$ of the left-hand side of \eqref{eq:classical_Whitney_2} is 
\begin{equation}
\sum_{\substack{0 \le k, l \le d \\ k+l = d}} \class{\extprod^{k} \CS} \cdot \class{\extprod^{l} \CQ}, 
\end{equation}
while that of the right-hand side of \eqref{eq:classical_Whitney_2} is $\class{\extprod^{d} \udBC^{2n}} = e_{d}(\e^{\pm \ve})$. 

\begin{proof}[Proof of Theorem~\ref{thm:Whitney_presentation}]
Apply Theorem~\ref{thm:Nakayama} in the following situation: 
\begin{itemize}
\item $G = \Sp{2n}$ and $L = I \setminus \{n\}$, 
\item $(m_{1}, \ldots, m_{s}) = (e_{1}(x), \ldots, e_{n}(x), e_{1}(y), \ldots, e_{n}(y))$, 
\item $\ti{P}_{d} = \displaystyle\sum_{\substack{k, l \ge 0 \\ k+l=d}} e_{k}(x) e_{l}(y) - e_{d}(\e^{\pm \ve}) + QF_{d-2}$ for $1 \le d \le 2n$, 
\item $k = 1$ and $Q_{1} = Q$, and 
\item $\Phi = \Psi_{\cl}$. 
\end{itemize}
Then we obtain the lift $\Psi = \ti{\Psi}_{\cl}: (R(T)\bra{Q})[x_{1}, \ldots, x_{n}, y_{1}, \ldots, y_{n}]^{\FS_{n}}/J_{\LG(n)} \xrightarrow{\sim} QK_{T}(\LG(n))$, as desired. 
\end{proof}

\end{document}